\documentclass[11pt]{amsart}
\usepackage{fullpage}
\usepackage{amsmath, amssymb, amsthm}
\usepackage{algorithm}
\usepackage {algpseudocode}
\usepackage{graphicx}
\usepackage{color}
\usepackage {bm}  
\usepackage[colorlinks=true,linkcolor=red,citecolor=blue,urlcolor=cyan]{hyperref}
\usepackage{verbatim}

\newtheorem{thm}{Theorem}[section]

\newtheorem{lem}[thm]{Lemma}

\theoremstyle{remark}

\usepackage[backend=bibtex,giveninits=true,sorting=nyt,style=alphabetic,natbib=true,maxcitenames=2,maxbibnames=10,url=false,doi=true,backref=false]{biblatex}
\renewbibmacro{in:}{\ifentrytype{article}{}{\printtext{\bibstring{in}\intitlepunct}}}

\makeatletter
\DeclareFontFamily{U}{tipa}{}
\DeclareFontShape{U}{tipa}{m}{n}{<->tipa10}{}
\newcommand{\arc@char}{{\usefont{U}{tipa}{m}{n}\symbol{62}}}%

\newcommand{\arc}[1]{\mathpalette\arc@arc{#1}}

\newcommand{\arc@arc}[2]{%
  \sbox0{$\m@th#1#2$}%
  \vbox{
    \hbox{\resizebox{\wd0}{\height}{\arc@char}}
    \nointerlineskip
    \box0
  }%
}
\makeatother

\usepackage{pgf}

\newcommand{\IR}{\mathbb{R}}

\newcommand{\bkappa}{\bm{\kappa}}
\newcommand{\bPhi}{\bm{\Phi}}
\newcommand{\bb}{\mathbf{b}}
\newcommand{\bd}{\mathbf{d}}
\newcommand{\bff}{\mathbf{f}}
\newcommand{\bg}{\mathbf{g}}
\newcommand{\bp}{\mathbf{p}}
\newcommand{\bq}{\mathbf{q}}
\newcommand{\br}{\mathbf{r}}

\newcommand{\bu}{\mathbf{u}}
\newcommand{\bv}{\mathbf{v}}

\newcommand{\bz}{\mathbf{z}}
\newcommand{\BA}{\mathbf{A}}
\newcommand{\BC}{\mathbf{C}}
\newcommand{\BF}{\mathbf{F}}

\newcommand{\BP}{\mathbf{P}}
\newcommand{\BQ}{\mathbf{Q}}

\newcommand{\CA}{\mathcal{A}}
\newcommand{\CC}{\mathcal{C}}
\newcommand{\CF}{\mathcal{F}}
\newcommand{\CG}{\mathcal{G}}

\newcommand{\CH}{\mathcal{H}}
\newcommand{\CI}{\mathcal{I}}
\newcommand{\CJ}{\mathcal{J}}
\newcommand{\CK}{\mathcal{K}}
\newcommand{\CL}{\mathcal{L}}

\newcommand{\CN}{\mathcal{N}}
\newcommand{\CP}{\mathcal{P}}
\newcommand{\CR}{\mathcal{R}}
\newcommand{\CT}{\mathcal{T}}

\begin{document}
\title{A two-stage method for reconstruction of parameters in diffusion equations}

\author{Xuesong Bai}
\author{Elena Cherkaev}
\address{Department of Mathematics, University of Utah, Salt Lake City, UT}
\email{\{bai,elena\}@math.utah.edu}
\author{Dong Wang}
\address{School of Science and Engineering, The Chinese University of Hong Kong, Shenzhen, Guangdong, China}
\email{wangdong@cuhk.edu.cn}


\date{\today}

\begin{abstract} 
Parameter reconstruction for diffusion equations has a wide range of applications. In this paper, we proposed a two-stage  scheme to efficiently solve conductivity reconstruction problems for steady-state diffusion equations with solution data measured inside the domain.  The first stage is based on total variation regularization of the log diffusivity and the split Bregman iteration method. In the second stage, we apply the K-means clustering for the reconstruction of ``blocky'' conductivity functions. The convergence of the scheme is theoretically proved and extensive numerical examples are shown to demonstrate the performance of the scheme.
\end{abstract}

\maketitle

\section{Introduction}
Let $\Omega\subset\mathbb{R}^d$ be a bounded open set with piece-wise smooth boundary. Consider the following boundary value problem for a steady-state diffusion equation:
\begin{subequations} \label{bvp1}
 \begin{equation}
  -\nabla\cdot(\kappa(x)\nabla u(x)) = f(x), \quad x\in\Omega,
  \end{equation}
   \begin{equation}
  u(x) = g_D(x), \quad x\in\partial_D\Omega,
  \end{equation}
     \begin{equation}
  \kappa(x)\nabla u(x) \cdot \mathbf{n}(x) = g_N(x), \quad x\in\partial_N\Omega. \end{equation}
\end{subequations}
where $\mathbf{n}$ is the normal vector of the boundary.

Equation (\ref{bvp1}) arise in many physical problems such as thermal conduction \cite{Le_2017}, electrical conductivity, dielectrics, magnetism \cite{milton_2002}, optical tomography \cite{Gryazin_1999} and quantitative photo-acoustic tomography \cite{Bal_2013}, ground water flow \cite{vogel1999sparse}, and tumor growth modeling \cite{Gholami2016}.  

The parameter of interest is the conductivity $\kappa(x)$, or equivalently the log conductivity $q(x)=\ln(\kappa(x))$, which is to be reconstructed from data $z$ of the solution $u$. We would like $q$ to be square-integrable on $\Omega$, that is, $q\in L^2(\Omega)$, because the fact that $L^2(\Omega)$ is a Hilbert space will be needed in the analysis of the reconstruction scheme in this paper. However, we would not assume $q$ to be differentiable or even continuous. We are particularly interested in reconstructing ``blocky'' conductivity functions with jump discontinuities. Such discontinuous conductivity functions are often encountered in real world problems, and, in general, they are more difficult to reconstruct than continuous conductivity functions.

In this paper, the data $z$ is given as full-field measurements of the physical quantity $u$ ({\it i.e.}; measurements of the spatial distributions of $u$ inside the domain $\Omega$ \cite{grediac2012full}). Different approaches have been previously developed to reconstruct the conductivity in equation (\ref{bvp1}) from full-field measurements. In \cite{vogel1999sparse, vogel2002computational}, Newton-type methods were applied to minimize a penalized least-squares fit-to-data functional to reconstruct the conductivity. In \cite{bellis2017reconstructing}, authors used an iterative deconvolution method. Stochastic methods such as a randomized misfit approach can also be used (See \cite{Le_2017}).

In this paper, we follow the penalized least-squares approach and use notations from \cite{vogel1999sparse} to illustrate the formulation of this inverse problem. The PDE (\ref{bvp1}) can be represented as 
 \begin{equation}
  \CA(q)u=f
 \end{equation}
where $\CA(q)=-\nabla\cdot(\kappa\nabla)=-\nabla\cdot(\exp(q)\nabla)$.
 
Given a solution $u$ to the PDE (\ref{bvp1}), we denote the observation of the solution by
 \begin{equation}
  u_{obs}=\CC u
 \end{equation}
where $\CC$ is the observation operator. Then, the state-to-observation map can be represented as
 \begin{equation}
  \CF(q)=u_{obs}(q)= \CC u(q)=\CC\CA(q)^{-1}f.
 \end{equation}
Then the actual measured data, which is inexact, is written as
 \begin{equation}
  \label{data} z=\CF(q)+\eta
 \end{equation}
where $\eta$ represents the noise in the data. 
 
The reconstruction problem is now formulated as follows. Given measured data $z$, solve the penalized least-squares minimization problem
 \begin{equation}
  \label{pls1}\min_q \CJ_{pls}(q)\quad\mbox{where}\quad \CJ_{pls}(q) = \alpha\CJ_{reg}(q)+\frac{1}{2}||\CF(q)-z||_2^2.
 \end{equation}
Here $\CJ_{reg}$ is the regularization or penalty functional. In this paper, $\CJ_{reg}$ is chosen to be the isotropic total variation (TV) functional, which we will simply refer to as ``TV functional'' later:
 \begin{equation}\label{tv1}
  \CJ_{reg}(q)=TV(q)=||\nabla q||_1=\int_\Omega|\nabla q|dx=\int_\Omega\sqrt{\sum_{j=1}^d\left(\dfrac{\partial q}{\partial x_j}\right)^2}dx.
 \end{equation}
 
 For ``blocky'' and discontinuous conductivity $\kappa(x)$, the isotropic TV functional regularizes the surface area of the interfaces between blocks of different the conductivity. When the dimension $d=2$, this surface area is reduced to arc length. For example, it is shown in \cite{vogel2002computational} that if
 \begin{equation}
  q(x) = \left\{
  \begin{array}{ll}
   q_0, & q(x)\in\Omega_0, \\
   q_1, & q(x)\in\Omega\setminus\Omega_0,
  \end{array}
 \right.
 \label{block}
 \end{equation}
 where $\Omega_0\subset\Omega$, $q_0>q_1$, and the boundary $\partial\Omega_0$ is $C^2$, then
 \begin{equation}
 TV(q)=(q_0-q_1)Area(\partial\Omega_0).
 \end{equation}

Alternatively, the anisotropic TV functional may also be used:
 \begin{equation}
  \CJ_{reg}(q)=\sum_{j=1}^d\int_\Omega\left|\dfrac{\partial q}{\partial x_j}\right|dx.
 \end{equation}

Previously \cite{vogel1999sparse}, quasi-Newton methods are used to solve the penalized least-squares problem (\ref{pls1}). However, due to the nondifferentiability of the TV functional (\ref{tv1}) near $\nabla q = \mathbf{0}$, an approximation of the TV functional has to be used in this case:
\[
 \int_\Omega\sqrt{|\nabla q|^2+\beta}dx,
\]
where $\beta>0$ is a small number. Also, since it requires inverting the approximate Hessian of the objective functional, as $\alpha\to 0$, the condition number of the approximate Hessian approaches infinity, making the computation inefficient \cite{nocedal2006numerical}. For $L^1$-regularized problems such as the case of TV regularization, the split Bregman method serves as a better alternative \cite{goldstein2009split}.
  
Another challenge that arises in the reconstruction of ``blocky'' and discontinuous conductivity $\kappa(x)$ is that solutions of the discretized reconstruction problem are actually continuous, which cannot recover the discontinuous nature of such $\kappa(x)$. In these cases, a K-means clustering and thresholding process is applied to the continuous solution to reconstruct the ``blocks'' of the conductivity. This stage is motivated by the two-stage method for segmentation in \cite{cai2013two,Cai_2017,Wu_2021} and the iterative convolution thresholding method in  \cite{Wang_2017, Wang_2019}.

\section{The reconstruction scheme}

The challenge in the penalized least-squares minimization problem \eqref{pls1} for conductivity reconstruction is that the above functional has both $L^1$ and $L^2$ portions, making the problem difficult to solve numerically. The split Bregman method for $L^1$-regularized problems \cite{goldstein2009split} aims exactly at ``decoupling'' the $L^1$ and $L^2$ portion of $L^1$-regularized problems, which makes efficient numerical solvers possible.

There are several ways to derive and prove convergence of the split Bregman iteration scheme. For example, in \cite{goldstein2009split} the iteration scheme is derived from a more general Bregman iterative regularization scheme \cite{osher2005iterative,yin2007bregman} for minimization problems of the form
\begin{equation}
 \label{grm}\min_q\{\CJ_{reg}(q)+\CI(q,z)\},
\end{equation}
where $\CJ_{reg}$ is a convex nonnegative regularization functional and the fit-to-data functional $\CI$ is convex and nonnegative with respect to $q$ for fixed $z$. However, due to nonlinearity of $\CF$, the least-squares fit-to-data functional
\begin{equation}
 \CI(q,z) = \frac12||\CF(q)-z||_2^2
\end{equation}
may not be convex in (\ref{pls1}). As a result of this, the split Bregman scheme for (\ref{pls1}) cannot be derived this way, and alternative derivation is needed.

In this section, we first present the details of the split Bregman scheme applied to the conductivity reconstruction problem (\ref{pls1}). Then we show that the split Bregman scheme for least-squares minimization problems can be derived from the alternating direction method of multipliers (ADMM) \cite{esser2009applications,nien2014convergence} with nonlinear equality constraints. Finally, we give a convergence result for the nonlinear equality-constrained ADMM, which then shows the convergence of the split Bregman scheme.

\subsection{Split Bregman scheme for conductivity reconstruction problems}

First, we write (\ref{pls1}) (\ref{tv1}) in the form of
 \begin{equation}
  \label{l1rp} \min_q\left\{\frac{\mu}{2}||\CF(q)-z||_2^2+||\nabla q||_1\right\},
 \end{equation}
 where $\mu=\dfrac{1}{\alpha}$. The purpose of rewriting the problem is to make the coefficient in front of the $L^1$ portion to be $1$.
 
The split Bregman iteration scheme in \cite{goldstein2009split} solves the following constrained minimization problem equivalent to (\ref{l1rp}) using one split variable $\bd$:
\begin{equation}
 \label{cmp1} \min_{q,\bd}\left\{\frac{\mu}{2}||\CF(q)-z||_2^2+||\bd||_1\right\}\quad\mbox{s.t.}\quad \bd=\nabla q. 
\end{equation}
 
The split Bregman algorithm for constrained minimization problem (\ref{cmp1}) is given in Algorithm \ref{alg_1}.

\begin{algorithm}[H]
    \caption{Split Bregman iteration scheme for conductivity reconstruction}
    \begin{algorithmic}[1]
     \Require $\CF(q)$, $z$, $\mu$, $\lambda$.
     \State Initialize: $k=0$; assign values to $q^0$, $\bd^0$, and $\bb^0$.
     \While{``not converge''}
     \State $q^{k+1}\leftarrow\arg\min_q\left\{\dfrac{\mu}{2}||\CF(q)-z||_2^2+\dfrac{\lambda}{2}||\bd^k-\nabla q-\bb^k||_2^2\right\}$
     \label{state1}
     \State $\bd^{k+1}\leftarrow\arg\min_\bd\left\{||\bd||_1+\dfrac{\lambda}{2}||\bd-\nabla q^{k+1}-\bb^k||_2^2\right\}$ \label{state2}
     \State $\bb^{k+1}\leftarrow\bb^k+(\nabla q^{k+1} -\bd^{k+1})$ \label{state3}
     \State $k\leftarrow k+1$
     \EndWhile
    \end{algorithmic}
    \label{alg_1}
\end{algorithm}

The speed of this split Bregman scheme is largely dependent on how fast we can solve the two minimization subproblems.
 
 \subsubsection{q-subproblem}
 
 Since the $L^1$ and $L^2$ portions of (\ref{l1rp}) are now decoupled, the $q$-subproblem no longer contains any $L^1$ component; hence it is now differentiable. The nonlinear conjugate gradient method or BFGS method can be applied. However, the nonlinear conjugate gradient method \cite{nocedal2006numerical} does not give good results, possibly due to the nonlinearity of the operator $\CF$ as well as the large number of variables. 
 
 In this paper, MATLAB package \textit{Hybrid Algorithm for Nonsmooth Optimization} (HANSO) by Overton \cite{lewis2013nonsmooth} is used. The algorithm is based on BFGS and gradient sampling.
 
 \subsubsection{d-subproblem}
 
 We can explicitly compute the components $d_j^{k+1}$, $j=1,\cdots,d$ of $\bd^{k+1}$ using a generalized shrinkage formula \cite{goldstein2009split}:
 \begin{equation}\label{shrinkage}
  d_j^{k+1}=\max\left\{s^k-\frac{\lambda}{2},0\right\}\dfrac{\frac{\partial q^k}{\partial x_j}+b_j^k}{s^k},
 \end{equation}
 where
 \begin{equation}
  s^k=|\nabla q^{k+1}+\bb^k|=\sqrt{\sum_{j=1}^d\left(\frac{\partial q_k}{\partial x_j}+b_j^k\right)^2}.
 \end{equation}
 Indeed, the Euler-Lagrange equation of the d-subproblem is
 \begin{equation}
  \frac{\bd}{|\bd|}+\lambda(\bd-(\nabla q^{k+1}+\bb^k))=\mathbf{0},
 \end{equation}
 which has an explicit solution
 \begin{equation}
  \bd = \max\left\{s^k-\frac{\lambda}{2},0\right\}\frac{\nabla q^{k+1}+\bb^k}{s^k}.
 \end{equation}

\subsection{Segmentation scheme for the reconstruction of ``blocky'' conductivity}

 To better reconstruct ``blocky'' and discontinuous conductivity, we introduce an extra segmentation scheme by thresholding the continuous reconstruction result $\kappa(x)=\exp(q(x))$ to get a discontinuous and more realistic reconstruction $\kappa_d(x)$. Here we apply the segmentation scheme to the original conductivity function $\kappa(x)$ and not the log conductivity function $q(x)$, because the log scaling may distort the segmentation.
 
 In this paper, we use MATLAB's \texttt{kmeans} to cluster the reconstructed discretized $\bkappa(p) = \{\kappa(p)\} = \{\exp(q(p))\}$ from the numerical implementation of the split Bregman scheme (Algorithm \ref{alg_1}). Here, the set $\CP = \{p\}$ is either the set of cell centers in the cell centered finite difference (CCFD) method or the set of meshpoints in the finite element method (FEM). We will further explain these discretization methods later in Section \ref{discretization}
 
 Suppose that we want to segment $\bkappa(p)$ into $K$ segments, $K\geqslant 2$. Following the procedure described in \cite{cai2013two}, we use \texttt{kmeans} to classify the discrete set $\{\kappa(p)\}$ into $K$ clusters $\CK_1,\CK_2,\cdots\CK_K$, with mean value of each cluster $\hat{\rho}_1\leqslant\hat{\rho}_2\leqslant\cdots\leqslant\hat{\rho}_K$, without loss of generality. Then define the $(K-1)$ thresholds as
 \begin{equation}
  \rho_i = \frac{\hat{\rho}_i+\hat{\rho}_{i+1}}{2},\quad i=1,\cdots,K-1.
 \end{equation}
 The $i$th phase of $\kappa(p)$, $1\leqslant i \leqslant K$, is given by
 \begin{equation}
  \CP_i = \{p\in\CP:\rho_{i-1}<\kappa(p)\leqslant\rho_i\},
 \end{equation}
 where $\rho_0$ is set to be $0$.
 
 Alternatively, the $i$th phase of $\bkappa(p)$ can be directly given by
 \begin{equation}
  \CP_i = \{p\in\CP:\kappa(p)\in\CK_i\}
 \end{equation}
 without the thresholding procedure.
 
 The discontinuous reconstruction $\bkappa_d(p) = \{\kappa_d(p)\}$ is then set to be
 \begin{equation}
 \kappa_d(p) = \sum_{i=1}^K \hat{\rho}_i\chi_{\CP_i}(p),  
 \end{equation}
 where $\chi_{\CP_i}$ is the characteristic function of the set $\CP_i$.
 
 Note that in our current segmentation scheme, the number $K$ of different segments has to be set manually in advance and depends on \textit{a priori} knowledge about the conductivity function $\kappa(x)$. It may be possible to further automate the segmentation scheme by incorporating algorithms that analyze the pattern of the reconstruction result $\bkappa(p)$ and determine the number $K$ of segments. This will be a topic for future research.

\section{Derivation and convergence of the nonlinear split Bregman scheme}
 
\subsection{The split Bregman method as a nonlinear-equality constrained ADMM method}

In this subsection, we will show how the split Bregman scheme applied to the conductivity reconstruction problem (\ref{l1rp}), which is a nonlinear least-squares minimization problem, can be derived as an ADMM scheme with nonlinear equality constraints. The discussion here is largely inspired by \cite{nien2014convergence} and lays the foundation of the convergence analysis that will be presented in the next subsection.

\subsubsection{ADMM method for conductivity reconstruction problems}

Similar to the split Bregman method, the ADMM method also solves a constrained minimization problem equivalent to (\ref{l1rp}), but it uses two split variables $s$ and $\bd$ instead of one: 
\begin{equation}
 \label{cmp2} \min_{s,\bd}\left\{\frac{\mu}{2}||s-z||^2_2+||\bd||_1\right\}\quad\mbox{s.t.}\quad s=\CF(q),\,\bd=\nabla q.
\end{equation}

The ADMM algorithm for constrained minimization problem (\ref{cmp2}) is given in Algorithm \ref{alg_2}. The only difference between Algorithm \ref{alg_2} and the ADMM algorithm in \cite{nien2014convergence} is that the constraint $s=\CF(q)$ is nonlinear.

\begin{algorithm}[H]
    \caption{ADMM iteration scheme for conductivity reconstruction}
    \begin{algorithmic}[1]
     \Require $\CF(q)$, $z$, $\mu$, $\rho$, $\lambda$.
     \State Initialize: $k=0$; assign values to $q^0$, $s^0$, $\bd^0$, $c^0$, and $\bb^0$.
     \While{``not converge''}
     \State $q^{k+1}\leftarrow\arg\min_q\left\{\dfrac{\rho}{2}||s^k-\CF(q)-c^k||_2^2+\dfrac{\lambda}{2}||\bd^k-\nabla q-\bb^k||_2^2\right\}$
     \label{statea1}
     \State $s^{k+1}\leftarrow\arg\min_s\left\{\dfrac{\mu}{2}||s-z||_2^2+\dfrac{\rho}{2}||s-\CF(q^{k+1})-c^k||_2^2\right\}$ \label{statea2}
     \State $\bd^{k+1}\leftarrow\arg\min_\bd\left\{||\bd||_1+\dfrac{\lambda}{2}||\bd-\nabla q^{k+1}-\bb^k||_2^2\right\}$ \label{statea3}
     \State $c^{k+1}\leftarrow c^k+(\CF(q^{k+1})-s^{k+1})$ \label{statea4}
     \State $\bb^{k+1}\leftarrow\bb^k+(\nabla q^{k+1} -\bd^{k+1})$ \label{statea5}
     \State $k\leftarrow k+1$
     \EndWhile
    \end{algorithmic}
    \label{alg_2}
\end{algorithm}

\subsubsection{Derivation of the split Bregman scheme from the ADMM scheme}

Now we wish to show that the split Bregman scheme (Algorithm \ref{alg_1}) can be derived from the ADMM scheme (Algorithm \ref{alg_2}). Note that the $s$-subproblem in Algorithm \ref{alg_2} has a closed-form solution
\begin{equation}
 \label{ssub} s^{k+1} = \frac{\mu}{\mu+\rho}z+\frac{\rho}{\mu+\rho}(\CF(q^{k+1})+c^k).
\end{equation}
Combine (\ref{ssub}) with the $c$-subproblem, we see that if we initialize Algorithm \ref{alg_2} so that 
\begin{equation}
 s^0-\frac{\rho}{\mu}c^0 = z.
\end{equation}
 then we have the following identity for each $k$:
\begin{equation}
 \label{scid} s^{k+1}-\frac{\rho}{\mu}c^{k+1}=z,
\end{equation}
Substitute (\ref{scid}) into Algorithm \ref{alg_2} and we get the simplified ADMM algorithm for conductivity reconstruction, presented in Algorithm \ref{alg_3}.

\begin{algorithm}[H]
    \caption{Simplified ADMM iteration scheme for conductivity reconstruction}
    \begin{algorithmic}[1]
     \Require $\CF(q)$, $z$, $\mu$, $\rho$, $\lambda$.
     \State Initialize: $k=0$; assign values to $q^0$, $s^0$, $\bd^0$, $c^0$, and $\bb^0$, so that $s^0-\dfrac{\rho}{\mu}c^0 = z$.
     \While{``not converge''}
     \State $q^{k+1}\leftarrow\arg\min_q\left\{\dfrac{\rho}{2}\left|\left|\left(1-\dfrac{\mu}{\rho}\right)s^k+\dfrac{\mu}{\rho}z-\CF(q)\right|\right|_2^2+\dfrac{\lambda}{2}||\bd^k-\nabla q-\bb^k||_2^2\right\}$
     \label{states1}
     \State $s^{k+1}\leftarrow\dfrac{\rho}{\mu+\rho}\CF(q^{k+1})+\dfrac{\mu}{\mu+\rho}s^k$ \label{states2}
     \State $\bd^{k+1}\leftarrow\arg\min_\bd\left\{||\bd||_1+\dfrac{\lambda}{2}||\bd-\nabla q^{k+1}-\bb^k||_2^2\right\}$ \label{states3}
     \State $\bb^{k+1}\leftarrow\bb^k+(\nabla q^{k+1} -\bd^{k+1})$ \label{states4}
     \State $k\leftarrow k+1$
     \EndWhile
    \end{algorithmic}
    \label{alg_3}
\end{algorithm}

If we set $\rho=\mu$ in Algorithm \ref{alg_3}, then the $q$-subproblem is independent of $s^k$, and the simplified ADMM scheme is exactly reduced to the split Bregman scheme (Algorithm \ref{alg_1}). Thus we can show the convergence of the split Bregman scheme by proving the convergence of the ADMM scheme (Algorithm \ref{alg_2}).

\subsection{Convergence results for nonlinear-equality constrained ADMM method}

In this section, we follow the proofs in \cite{wang2017nonconvex} and give convergence results for ADMM schemes for nonlinear-equality constrained minimization problem with two split variables of the form (\ref{cmp2}), which then proves the convergence of the split Bregman scheme (Algorithm \ref{alg_1}).

\subsubsection{General ADMM scheme for nonlinear-equality constrained minimization problems with two split variables}

Consider minimization problems of the form
\begin{equation}
 \label{cmpg} \min_{s,\bd}\{\CH(s)+\CG(\bd)\}\quad\mbox{s.t.}\quad s=\CF(q), \bd=\bPhi(q),
\end{equation}
where $\CH(s)$ and $\CG(\bd)$ are convex functionals and the constraint functions $\CF(q)$ and $\bPhi(q)$ can be nonlinear.

The augmented Lagrangian for problem (\ref{cmpg}) \cite{esser2009applications,wang2017nonconvex} is
\begin{equation}
 \label{alag} L_{\lambda,\rho}(s,\bd,q,\bb,c) = \CH(s) + \CG(\bd) + \langle c,\CF(q)-s \rangle + \langle \bb, \bPhi(q)-\bd \rangle + \frac{\rho}{2}||s-\CF(q)||_2^2 + \frac{\lambda}{2}||\bd-\bPhi(q)||_2^2.
\end{equation}
The ADMM iteration scheme for (\ref{cmpg}) is then given in Algorithm \ref{alg_4}.

\begin{algorithm}[H]
    \caption{ADMM iteration scheme for nonlinear equality-constrained minimization problem with two split variables}
    \begin{algorithmic}[1]
     \Require $\CH(\cdot)$, $\CG(\cdot)$, $\CF(\cdot)$, $\bPhi(\cdot)$; $\rho$, $\lambda$.
     \State Initialize: $k=0$; assign values to $q^0$, $s^0$, $\bd^0$, $c^0$, and $\bb^0$.
     \While{``not converge''}
     \State $q^{k+1}\leftarrow\arg\min_q L_{\lambda,\rho}(s^k,\bd^k,q,\bb^k,c^k)$
     \label{stateg1}
     \State $s^{k+1}\leftarrow\arg\min_s L_{\lambda,\rho}(s,\bd^k,q^{k+1},\bb^k,c^k)$ \label{stateg2}
     \State $\bd^{k+1}\leftarrow\arg\min_\bd L_{\lambda,\rho}(s^{k+1},\bd,q^{k+1},\bb^k,c^k)$ \label{stateg3}
     \State $c^{k+1}\leftarrow c^k+\rho(\CF(q^{k+1})-s^{k+1})$ \label{stateg4}
     \State $\bb^{k+1}\leftarrow\bb^k+\lambda(\bPhi(q^{k+1}) -\bd^{k+1})$ \label{stateg5}
     \State $k\leftarrow k+1$
     \EndWhile
    \end{algorithmic}
    \label{alg_4}
\end{algorithm}

It is easy to check that if $\CH(s)=\frac{\mu}{2}||s-z||_2^2$, $\CG(\bd)=||\bd||_1$, and $\bPhi(q)=\nabla q$, then the iteration scheme given by Algorithm \ref{alg_4} is indeed equivalent to Algorithm \ref{alg_2}. The only apparent differences are the updates of $c^k$ (\ref{stateg4}) and $\bb^k$ (\ref{stateg5}). However, these differences do not affect the equivalence of the two algorithms. Indeed, we could slightly modify the augmented Lagrangian (\ref{alag}):
\begin{equation}\label{alag1}
 L_{\lambda,\rho}(s,\bd,q,\bb,c) = \CH(s) + \CG(\bd) + \rho\langle c,\CF(q)-s \rangle + \lambda\langle \bb, \bPhi(q)-\bd \rangle + \frac{\rho}{2}||s-\CF(q)||_2^2 + \frac{\lambda}{2}||\bd-\bPhi(q)||_2^2.
\end{equation}
(\ref{alag}) and (\ref{alag1}) are equivalent, and we can apply the the updates of $c^k$ (\ref{statea4}) and $\bb^k$ (\ref{statea5}) in Algorithm \ref{alg_2} to (\ref{alag1}).

\subsubsection{Theoretical analysis of convergence}

We begin the analysis by first specifying some notations. The (unaugmented) Lagrangian for problem (\ref{cmpg}) is given by
\begin{equation}
 \label{lag} L(s,\bd,q,\bb,c) = \CH(s) + \CG(\bd) + \langle c,\CF(q)-s \rangle + \langle \bb, \bPhi(q)-\bd \rangle.
\end{equation}
Then a saddle point of (\ref{lag}) is a point $(s^*, \bd^*, q^*, \bb^*, c^*)$ that satisfies
\begin{equation}
 \label{saddle} L(s^*,\bd^*,q^*,\bb,c)\leqslant L(s^*,\bd^*,q^*,\bb^*,c^*) \leqslant L(s,\bd,q,\bb^*,c^*)\quad\mbox{for any}\quad (s,\bd,q,\bb,c).
\end{equation}
Finding optimal solutions of (\ref{cmpg}) is equivalent to finding a saddle point of (\ref{lag}), which is then equivalent to the Kuhn-Tucker optimality conditions  \cite{rockafellar1970convex}.

Now we introduce some more notations. Let
\begin{equation}
 \CJ^k=\CH(s^k)+\CG(\bd^k)
\end{equation}
be the value of the objective functional in (\ref{cmpg}) at the $k$-th iteration, and let
\begin{equation}
 \CR_s^k=\CF(q^k)-s^k\quad\mbox{and}\quad\CR_\bd^k=\bPhi(q)^k-\bd^k
\end{equation}
be the residuals with respect to $s$ and $\bd$ at the $k$-th iteration, respectively. Corresponding to the above notations, we have
\begin{equation}
 \CJ^*=\CH(s^*)+\CG(\bd^*),\quad\mbox{and}\quad \CR_s^*=\CF(q^*)-s^*=0,\quad \CR_d^*=\bPhi(q^*)-\bd^*=\mathbf{0},
\end{equation}
where $\CJ^*$ is the optimal value of the objective functional in (\ref{cmpg}) at the optimal solution $(s^*, \bd^*, q^*)$, and the last two identities come from the Kuhn-Tucker optimality conditions.

The following lemma gives an estimation of the differences $\CJ^*-\CJ^{k+1}$ and $\CJ^{k+1}-\CJ^*$.

\begin{lem} \label{estim}
 Suppose that there exists an optimal solution $(s^*,\bd^*,q^*,\bb^*,c^*)$ of problem (\ref{cmpg}). Then
 \begin{equation}
  \label{Jestim1} \CJ^*-\CJ^{k+1}\leqslant \langle c^*, \CR_s^{k+1} \rangle + \langle \bb^*, \CR_\bd^{k+1} \rangle,
 \end{equation}
 and
 \begin{align}
  \label{Jestim2}\CJ^{k+1}-\CJ^* \leqslant& - \rho \langle s^{k+1}-s^k, \CF(q^{k+1})-\CF(q^*) \rangle \\
  & - \lambda \langle \bd^{k+1}-\bd^k, \bPhi(q^{k+1})-\bPhi(q^*) \rangle - \langle c^{k+1}, \CR_s^{k+1} \rangle - \langle \bb^{k+1}, \CR_\bd^{k+1} \rangle. \nonumber
 \end{align}
\end{lem}

From Lemma \ref{estim}, we get the following theorem about convergence properties of Algorithm \ref{alg_4}.

\begin{thm} \label{convergence}
 Suppose that there exists an optimal solution $(s^*,\bd^*,q^*,\bb^*,c^*)$ of problem (\ref{cmpg}). Then Algorithm \ref{alg_4} satisfies the following convergence properties:
 \begin{equation}
  \CR_s^k\to 0\quad\mbox{and}\quad\CR_\bd^k\to\mathbf{0}\quad\mbox{as}\quad k\to\infty\quad\mbox{strongly},
 \end{equation}
 and 
 \begin{equation}
  \CJ^k\to\CJ^*\quad\mbox{as}\quad k\to\infty.
 \end{equation}
\end{thm}

Theorem \ref{convergence} guarantees that the residuals with respect to $s$ and $\bd$ will go to zero, and that the value of objective functional in (\ref{cmpg}) will converge to an optimal value. Note, however, that Theorem \ref{convergence} does not guarantee the convergence of $(s^k,\bd^k,q^k,\bb^k,c^k)$ to the specified optimal solution $(s^*,\bd^*,q^*,\bb^*,c^*)$.

\section{Proof of convergence results}

In this section we give a proof for Lemma \ref{estim} and Theorem \ref{convergence}. The proof follows the methods presented in \cite{wang2017nonconvex}.

\subsection{Notations and lemmas}

Let $X$ and $Y$ be real Banach spaces, and let $X^*$ and $Y^*$ denote the corresponding dual spaces. For a convex functional $\CI_c:X\to\IR$, the (regular) subdifferential at a point $q$ is defined as 
 \begin{equation}
  \partial\CI_c(q)=\{q^*\in X^*:\CI_c(p)\geqslant\CI_c(q)+\langle q^*, p-q\rangle,\,\forall p\in X\},
 \end{equation}
 where $\langle\cdot,\cdot\rangle$ denotes the standard duality product.

For a nonconvex functional $\CI_{nc}:X\to\IR$, the Fr\'{e}chet subdifferential \cite{kruger2003frechet}, defined below, acts as a generalization of the regular subdifferential:
\begin{equation}
 \partial\CI_{nc}(q)=\left\{q^*\in X^*:\liminf_{p\to q}\frac{\CI_{nc}(s)-\CI_{nc}(q)-\langle q^*, p-q\rangle}{||p-q||}\geqslant0\right\}.
\end{equation}
It is easy to show that if $\CI_{nc}$ is Fre\'{e}chet differentiable, then $\partial\CI_{nc}(q)=\{\CI_{nc}'(q)\}$, where $\CI_{nc}'(q)$ is the Fre\'{e}chet derivative of $\CI_{nc}$. We will use the same notation $\partial\CI(q)$ to denote either the regular or the Fr\'{e}chet subdifferential of a general functional $\CI:X\to\IR$.

We will need the following calculus results for Fr\'{e}chet subdifferentials. The first one is a sum rule, which follows directly from the definitions.
\begin{lem}\label{sum}
 Let $\CI_1$, $\CI_2$: $X\to\IR$ be subdifferentiable at $q$. Then $\CI_1+\CI_2$ is subdifferentiable at $q$ and
 \[
  \partial \CI_1(q) + \partial \CI_2(q) \subset \partial (\CI_1+\CI_2) (q).
 \]
\end{lem}
The next result is a chain rule. We will need some more notations before introducing this result. Let $F:X\to Y$ be a function between the two Banach spaces. For any $s^*\in Y^*$, we can define a scalar function $\langle s^*, F\rangle:X\to\IR$ by
\[
 \langle s^*,F \rangle (q) = \langle s^*,F(q) \rangle.
\]
Now let $\CJ:Y\to\IR$ be a functional on $Y$ that is subdifferentiable at $s=F(q)$. Consider the composition $\CI(q) = \CJ(F(q))$. We have the following lemma \cite{kruger2003frechet,mordukhovich2006frechet}.
\begin{lem}\label{chain}
 Let $\CJ$ be Fr\'{e}chet differentiable at $s=F(q)$. Then
 \[
  \partial \CI(q) = \partial \langle \CJ'(s), F \rangle (q).
 \]
\end{lem}
In this paper, we take $X=Y=L^2(\Omega)$. In this case, $\langle \cdot,\cdot \rangle$ is the $L^2$ - inner product.

\subsection{Proof of Lemma \ref{estim}}

We begin by proving inequality (\ref{Jestim1}). This can be done by simply plugging the point $(s^{k+1},\bd^{k+1},q^{k+1},\bb^{k+1},c^{k+1})$ in the saddle point condition (\ref{saddle}):
\begin{equation} \label{saddlekp1}
 L(s^*,\bd^*,q^*,\bb^*,c^*) - L(s^{k+1},\bd^{k+1},q^{k+1},\bb^*,c^*) \leqslant 0.
\end{equation}
From (\ref{saddlekp1}) we directly obtain
\begin{equation*}
  \CJ^*-\CJ^{k+1}\leqslant \langle c^*, \CR_s^{k+1} \rangle + \langle \bb^*, \CR_\bd^{k+1} \rangle
 \end{equation*}
and inequality (\ref{Jestim1}) is proved.

To prove inequality (\ref{Jestim2}), we need the following series of inequalities. First, observe from line \ref{stateg1} of Algorithm \ref{alg_4} that
\begin{equation}
 \label{qopt1} q^{k+1} = \arg\min_q L_{\lambda,\rho}(s^k,\bd^k,q,\bb^k,c^k).
\end{equation}
From the optimality condition of equation (\ref{qopt1}), Lemma \ref{chain}, and Lemma \ref{sum}, we have that
\begin{equation} \label{qopt2}  \arraycolsep=1.4pt \def\arraystretch{2}
\begin{array}{lll}
 0 &\in & \partial\langle c^k,\CF(\cdot)-s^k \rangle (q^{k+1}) + \partial\langle \bb^k, \bPhi(\cdot)-\bd^k \rangle (q^{k+1}) + \partial\rho \langle \CF(q^{k+1})-s^k, \CF(\cdot)-s^k \rangle (q^{k+1}) \\
  & & + \partial\lambda \langle \bPhi(q^{k+1})-\bd^k, \bPhi(\cdot)-\bd^k \rangle (q^{k+1}) \\
 & \subset & \partial\langle c^k+\rho(\CF(q^{k+1})-s^k),\CF(\cdot)-s^k \rangle (q^{k+1}) + \partial\langle \bb^k+\lambda(\bPhi(q^{k+1})-\bd^k), \bPhi(\cdot)-\bd^k \rangle (q^{k+1}).
\end{array}
\end{equation}
Rewriting line \ref{stateg4} and \ref{stateg5} of Algorithm \ref{alg_4} we have
\begin{equation}
 \label{ckbk} c^k = c^{k+1}-\rho(\CF(q^{k+1})-s^{k+1})\quad\mbox{and}\quad \bb^k = \bb^{k+1}-\lambda(\bPhi(q^{k+1})-\bd^{k+1}).
\end{equation}
Substitute (\ref{ckbk}) into (\ref{qopt2}) we obtain
\begin{equation}
 \label{qopt3} 0\in \partial\langle c^{k+1}+\rho(s^{k+1}-s^k),\CF(\cdot)-s^k \rangle (q^{k+1}) + \partial\langle \bb^{k+1}+\lambda(\bd^{k+1}-\bd^k), \bPhi(\cdot)-\bd^k \rangle (q^{k+1}).
\end{equation}
This implies that
\begin{equation}
 \label{qopt4} q^{k+1} = \arg\min_q \{\langle c^{k+1}+\rho(s^{k+1}-s^k),\CF(q)-s^k \rangle + \partial\langle \bb^{k+1}+\lambda(\bd^{k+1}-\bd^k), \bPhi(q)-\bd^k \rangle\}.
\end{equation}
It follows from (\ref{qopt4}) that
\begin{equation}\label{qopt5} \arraycolsep=1.4pt \def\arraystretch{2.2} 
 \begin{array}{ll}
    & \langle c^{k+1},\CF(q^{k+1})-\CF(q^*) \rangle + \langle \bb^{k+1}, \bPhi(q^{k+1})-\bPhi(q^*) \rangle \\
  + & \rho \langle s^{k+1}-s^k,\CF(q^{k+1})-\CF(q^*) \rangle + \lambda \langle \bd^{k+1}-\bd^k,\bPhi(q^{k+1})-\bPhi(q^*) \rangle \leqslant 0.
 \end{array}
\end{equation}

Next, observe from line \ref{stateg2} of Algorithm \ref{alg_4} that
\begin{equation}
 \label{sopt1} s^{k+1} = \arg\min_s L_{\lambda,\rho}(s,\bd^k,q^{k+1},\bb^k,c^k).
\end{equation}
From the optimality condition of (\ref{sopt1}), we obtain that
\begin{equation}
 \label{sopt2} 0\in \partial\CH(s^{k+1}) - c^k - \rho (\CF(q^{k+1})-s^{k+1}).
\end{equation}
Substituting (\ref{ckbk}) into (\ref{sopt2}) gives
\begin{equation}
 \label{sopt3} 0\in \partial\CH(s^{k+1}) - c^{k+1}.
\end{equation}
This implies that
\begin{equation}
 \label{sopt4} s^{k+1} = \arg\min_s\{\CH(s) + \langle c^{k+1}, \CF(q^{k+1})-s \rangle \}.
\end{equation}
It follows from (\ref{sopt4}) that
\begin{equation*}
 \CH(s^{k+1}) - \CH(s^*) + \langle c^{k+1},\CF(q^{k+1})-s^{k+1} \rangle - \langle c^{k+1},\CF(q^{k+1})-s^* \rangle \leqslant 0.
\end{equation*}
That is,
\begin{equation}
 \label{sopt5} \CH(s^{k+1}) - \CH(s^*) + \langle c^{k+1}, \CR_s^{k+1} \rangle - \langle c^{k+1},\CF(q^{k+1})-s^* \rangle \leqslant 0.
\end{equation}
Analogously, from line \ref{stateg3} of Algorithm \ref{alg_4} we obtain
\begin{equation}
 \label{dopt5} \CG(\bd^{k+1}) - \CG(\bd^*) + \langle \bb^{k+1},\CR_\bd^{k+1} \rangle - \langle \bb^{k+1},\bPhi(q^{k+1})-\bd^* \rangle \leqslant 0.
\end{equation}
Adding inequalities (\ref{qopt5}), (\ref{sopt5}), (\ref{dopt5}) together and using the facts that $\CF(q^*)-s^*=0$ and $\bPhi(q^*)-\bd^*=\mathbf{0}$, we get
\begin{equation*} \arraycolsep=1.4pt \def\arraystretch{2.2} 
 \begin{array}{ll} 
    & \CJ^{k+1}-\CJ^* + \rho \langle s^{k+1}-s^k, \CF(q^{k+1})-\CF(q^*) \rangle + \lambda \langle \bd^{k+1}-\bd^k, \bPhi(q^{k+1})-\bPhi(q^*) \rangle \\
    & + \langle c^{k+1}, \CR_s^{k+1} \rangle + \langle \bb^{k+1}, \CR_\bd^{k+1} \rangle \leqslant 0,
 \end{array}
\end{equation*}
and inequality (\ref{Jestim2}) is proved.

\subsection{Proof of Theorem \ref{convergence}}  

Adding (\ref{Jestim1}) and (\ref{Jestim2}) together and then multiplying both sides by 2, we get
\begin{equation}\label{initineq} \arraycolsep=1.4pt \def\arraystretch{1.5} 
 \begin{array}{ll}
    & 2\langle c^{k+1}-c^*,\CR_s^{k+1} \rangle + 2\rho\langle s^{k+1}-s^k,\CF(q^{k+1})-\CF(q^*) \rangle \\
  + & 2\langle \bb^{k+1}-\bb^*,\CR_\bd^{k+1} \rangle + 2\lambda\langle \bd^{k+1}-\bd^k,\bPhi(q^{k+1})-\bPhi(q^*) \rangle \leqslant 0.
 \end{array}
\end{equation}
We rewrite the terms in (\ref{initineq}). First
\begin{equation} \label{lyas1} \arraycolsep=1.4pt \def\arraystretch{1.5} 
 \begin{array}{ll}
    & 2\langle c^{k+1}-c^*,\CR_s^{k+1} \rangle \\
  = & 2\langle c^k+\rho\CR_s^{k+1}-c^*,\CR_s^{k+1} \rangle \quad\quad (c^{k+1}=c^k+\rho\CR_s^{k+1}) \\
  = & 2\langle c^k-c^*, \CR_s^{k+1} \rangle + 2\rho||\CR_s^{k+1}||_2^2 \\
  = & \frac{2}{\rho}\langle c^k-c^*,c^{k+1}-c^k \rangle + \frac{1}{\rho}||c^{k+1}-c^k||_2^2 + \rho||\CR_s^{k+1}||_2^2 \quad\quad \left(\CR_s^{k+1}=\frac{c^{k+1}-c^k}{\rho}\right) \\
  = & \frac{2}{\rho}\langle c^k-c^*,(c^{k+1}-c^*)-(c^k-c^*) \rangle + \frac{1}{\rho}||(c^{k+1}-c^*)-(c^k-c^*)||_2^2 + \rho||\CR_s^{k+1}||_2^2 \\
  = & \frac{1}{\rho}(||c^{k+1}-c^*||_2^2 - ||c^k-c^*||_2^2) + \rho||\CR_s^{k+1}||_2^2
 \end{array}
\end{equation}
Now consider the term $2\rho\langle s^{k+1}-s^k,\CF(q^{k+1})-\CF(q^*) \rangle$ together with the term $\rho||\CR_s^{k+1}||_2^2$ from (\ref{lyas1}):
\begin{equation} \label{lyas2} \arraycolsep=1.4pt \def\arraystretch{1.5} 
 \begin{array}{ll}
    & 2\rho\langle s^{k+1}-s^k,\CF(q^{k+1})-\CF(q^*) \rangle + \rho||\CR_s^{k+1}||_2^2 \\
  = & 2\rho\langle s^{k+1}-s^k,s^{k+1}+\CR_s^{k+1}-s^* \rangle + \rho||\CR_s^{k+1}||_2^2 \\
  = & 2\rho\langle s^{k+1}-s^k,s^{k+1}-s^* \rangle + 2\rho\langle s^{k+1}-s^k,\CR_s^{k+1} \rangle + \rho||\CR_s^{k+1}||_2^2 \\
  = & 2\rho\langle s^{k+1}-s^k,s^{k+1}-s^k+s^k-s^* \rangle + 2\rho\langle s^{k+1}-s^k,\CR_s^{k+1} \rangle + \rho||\CR_s^{k+1}||_2^2 \\
  = & 2\rho||s^{k+1}-s^k||_2^2 + 2\rho\langle s^{k+1}-s^k,s^k-s^* \rangle + 2\rho\langle s^{k+1}-s^k,\CR_s^{k+1} \rangle + \rho||\CR_s^{k+1}||_2^2 \\
  = & (\rho||s^{k+1}-s^k||_2^2 + 2\rho\langle s^{k+1}-s^k,\CR_s^{k+1} \rangle + \rho||\CR_s^{k+1}||_2^2) + 2\rho\langle s^{k+1}-s^k,s^k-s^* \rangle \\
  + & \rho||s^{k+1}-s^k||_2^2 \\
  = & \rho||s^{k+1}-s^k+\CR_s^{k+1}||_2^2 + 2\rho\langle (s^{k+1}-s^*)-(s^k-s^*),s^k-s^* \rangle \\
  + & \rho||(s^{k+1}-s^*)-(s^k-s^*)||_2^2 \\
  = & \rho||s^{k+1}-s^k+\CR_s^{k+1}||_2^2 + \rho||s^{k+1}-s^*||_2^2 - \rho||s^k-s^*||_2^2 \\
 \end{array}
\end{equation}

From (\ref{lyas1}) and (\ref{lyas2}) we get that
\begin{equation}\label{lyas3}\arraycolsep=1.4pt \def\arraystretch{1.5} 
 \begin{array}{ll}
    & 2\langle c^{k+1}-c^*,\CR_s^{k+1} \rangle + 2\rho\langle s^{k+1}-s^k,\CF(q^{k+1})-\CF(q^*) \rangle \\
  = & \rho||s^{k+1}-s^k+\CR_s^{k+1}||_2^2 + \rho||s^{k+1}-s^*||_2^2 - \rho||s^k-s^*||_2^2 + \frac{1}{\rho}(||c^{k+1}-c^*||_2^2 - ||c^k-c^*||_2^2)
 \end{array}
\end{equation}
Analogously, for the other two terms in (\ref{initineq}) we have
\begin{equation}\label{lyad3}\arraycolsep=1.4pt \def\arraystretch{1.5} 
 \begin{array}{ll}
    & 2\langle \bb^{k+1}-\bb^*,\CR_\bd^{k+1} \rangle + 2\lambda\langle \bd^{k+1}-\bd^k,\bPhi(q^{k+1})-\bPhi(q^*) \rangle \\
  = & \lambda||\bd^{k+1}-\bd^k+\CR_\bd^{k+1}||_2^2 + \lambda||\bd^{k+1}-\bd^*||_2^2 - \lambda||\bd^k-\bd^*||_2^2 \\
  + & \frac{1}{\lambda}(||\bb^{k+1}-\bb^*||_2^2 - ||\bb^k-\bb^*||_2^2)
 \end{array}
\end{equation}

Let the Lyapunov function $V^k$ be defined as
\begin{equation}
 V^k = \rho||s^k-s^*||_2^2 + \lambda||\bd^k-\bd^*||_2 + \frac{1}{\rho}||c^k-c^*||_2^2 + \frac{1}{\lambda}||\bb^k-\bb^*||_2^2.
\end{equation}
Then from (\ref{initineq}), (\ref{lyas3}), and (\ref{lyad3}) we obtain
\begin{equation}\label{lyadecay1}
 V^{k+1}\leqslant V^k - \rho||s^{k+1}-s^k+\CR_s^{k+1}||_2^2 - \lambda||\bd^{k+1}-\bd^k+\CR_\bd^{k+1}||_2^2
\end{equation}
Recall from (\ref{sopt3}) that
\begin{equation}
 s^{k+1} = \arg\min_s\{\CH(s) + \langle c^{k+1}, -s \rangle \} \quad\mbox{and}\quad s^{k} = \arg\min_s\{\CH(s) + \langle c^{k}, -s \rangle \},
\end{equation}
so
\begin{equation} \label{hs1}
 \CH(s^{k+1}) + \langle c^{k+1}, -s^{k+1} \rangle \leqslant \CH(s^k) + \langle c^{k+1}, -s^k \rangle 
\end{equation}
and 
\begin{equation} \label{hs2}
 \CH(s^k) + \langle c^k, -s^k \rangle \leqslant \CH(s^{k+1}) + \langle c^k, -s^{k+1} \rangle. 
\end{equation}
Adding (\ref{hs1}) and (\ref{hs2}) together yields
\begin{equation}\label{cs1}
 \langle c^k-c^{k+1},s^{k+1}-s^k \rangle \leqslant 0.
\end{equation}
Substitute $c^k-c^{k+1}=-\rho\CR_s^{k+1}$ into (\ref{cs1}) we get
\begin{equation}
 \label{rs1}\langle \CR_s^{k+1},s^{k+1}-s^k \rangle \geqslant 0.
\end{equation}
Thus
\begin{equation}
 \label{sineq} - \rho||s^{k+1}-s^k+\CR_s^{k+1}||_2^2 \leqslant - \rho||s^{k+1}-s^k||_2^2 - \rho||\CR_s^{k+1}||_2^2.
\end{equation}
Analogously
\begin{equation}
 \label{dineq} - \rho||\bd^{k+1}-\bd^k+\CR_\bd^{k+1}||_2^2 \leqslant - \rho||\bd^{k+1}-\bd^k||_2^2 - \rho||\CR_\bd^{k+1}||_2^2.
\end{equation}
From inequalities (\ref{lyadecay1}), (\ref{sineq}), (\ref{dineq}) we get
\begin{equation}
 \label{lyadecay2} V^{k+1} \leqslant V^k - \rho||s^{k+1}-s^k||_2^2 - \rho||\CR_s^{k+1}||_2^2 - \rho||\bd^{k+1}-\bd^k||_2^2 - \rho||\CR_\bd^{k+1}||_2^2.
\end{equation}
Taking the infinite sum of (\ref{lyadecay2}) starting from $k=0$, we have
\begin{equation}\label{infsum}
 \rho\left(\sum_{k=0}^\infty||s^{k+1}-s^k||_2^2 + \sum_{k=0}^\infty||\CR_s^{k+1}||_2^2\right) + \lambda\left(\sum_{k=0}^\infty||\bd^{k+1}-\bd^k||_2^2 + \sum_{k=0}^\infty||\CR_\bd^{k+1}||_2^2\right) \leqslant V^0.
\end{equation}
This implies that
\begin{equation}
 \label{resconv} \lim_{k\to\infty} ||\CR_s^k||_2 = 0 \quad\mbox{and}\quad \lim_{k\to\infty} ||\CR_\bd^k||_2 = 0,
\end{equation}
hence
\begin{equation}
  \CR_s^k\to 0\quad\mbox{and}\quad\CR_\bd^k\to\mathbf{0}\quad\mbox{strongly}\quad \mbox{as}\quad k\to\infty.
\end{equation}
From the convergence result (\ref{resconv}) and the estimations (\ref{Jestim1}) and (\ref{Jestim2}) we get
\begin{equation}
 \lim_{k\to\infty} \CJ^k = \CJ^*.
\end{equation}

\section{Numerical methods}

To implement the split Bregman scheme and solve the conductivity reconstruction problem, we need to discretize the penalized least-squares problem (\ref{pls1})
\begin{equation}
 \min_{\bq\in\IR^N}\frac{1}{2}\|\BF(\bq)-\bz\|_2^2+\alpha J_{reg}(\bq),
\end{equation}
where 
\begin{equation}
 \BF(\bq) = \BC\bu(\bq) = \BC\BA(\bq)^{-1}\bff,
\end{equation}
$\BA(\bq)$ is an $N\times N$ stiffness matrix, and $\BC$ is an $N_c\times N$ matrix. The split Bregman iteration scheme (Algorithm \ref{alg_1}) is then discretized correspondingly.

\subsection{Discretization of the domain and solution to the diffusion equation}\label{discretization}

In this paper, we focus on the cases where $\Omega\subset\mathbb{R}^2$.

For the cases where $\Omega$ is a square or rectangular domain, cell-centered finite difference discretization (CCFD) \cite{arbogast1997mixed} is used to solve the boundary value problem (\ref{bvp1}). For other cases where the geometry of the domain $\Omega$ is more complicated, finite element method (FEM) with linear elements \cite{atkinson2005theoretical} \cite{strang2007computational} is needed.

A simple mesh generator {\tt distmesh2d} in MATLAB \cite{persson2004simple} is used along with FEM to triangulate the domain $\Omega$ into a triangular mesh $\CT=\{T\}$ with mesh points (or nodes) $\CP=\{p\}$. The geometry of the domain is represented by signed distance functions $dist(x,y)$, where $|dist(x,y)|$ is the distance from the point $(x,y)$ to the boundary $\partial\Omega$, and
\[
dist(x,y) = \left\{
   \begin{array}{ll}
    >0, & (x,y)\not\in\Omega, \\
    <0, & (x,y)\in\Omega.
   \end{array}
   \right.
\]
If the boundary $\partial \Omega$ of the domain is given piece-wise by explicit functions $y = h_i(x)$, then the corresponding signed distance function $dist_i(x,y)$ can be approximated by the following straightforward ``discretize-and-search'' algorithm:

\begin{algorithm}[H]
 \caption{Approximation of the signed distance function $dist_i(x,y)$}
 \begin{algorithmic}[1]
  \Require piece-wise boundary information $y = h_i(x)$ with domain $[X_1,X_2]$; $N$ as the number of nodes.
  \State Discretize the interval $[X_1,X_2]$ uniformly into $N$ nodes $\{x^j\}_{j=1}^N$.
  \State Calculate the $N$ distances $\{dist^j\}_{j=1}^N$ between $(x,y)$ and $(x^j,h_i(x^j))$. \label{discdist}
  \State 
  \[
   dist_i(x,y) \leftarrow \left\{
   \begin{array}{ll}
    \min_{1\leqslant j\leqslant N} dist^j, & (x,y)\not\in\Omega, \\
    -\min_{1\leqslant j\leqslant N} dist^j, & (x,y)\in\Omega.
   \end{array}
   \right.
  \]
 \end{algorithmic}
 \label{alg_5}
\end{algorithm}

The signed distance function $dist(x,y)$ for the whole domain $\Omega$ can then be constructed by using the codes \texttt{dunion}, \texttt{ddiff}, and \texttt{dintersect} to combine the geometries defined by $d_i(x,y)$'s. Despite its approximation nature of Algorithm \ref{alg_5}, it is shown that when $N$ is large enough, the signed distance functions it generates are accurate enough to work with \texttt{distmesh2d}. 

If $\partial\Omega$ is given by a set of discrete points instead, then we can either directly apply these points to step \ref{discdist} of Algorithm \ref{alg_5}, or interpolate these points to get the functions $y = h_i(x)$.

A modified version of the MATLAB finite element method code {\tt femcode2} is then applied to the triangular mesh.

\subsection{Calculation of gradient using adjoint method}

The HANSO algorithm requires a calculation of the gradient with respect to $\bq$ of the discretized functional in the subproblem \ref{state1}, which is done in two parts.
 
 The least-square gradient $\bg_{ls}(\bq)$ of the discrete fit-to-data functional $\dfrac{1}{2}||\BF(\bq)-\bz||_2^2$ is calculated by \cite{vogel2002computational}
 \begin{equation}
  \label{ls}[\bg_{ls}(\bq)]_j= \bv(\bq)^t\left[\frac{\partial \BA(\bq)}{\partial q_j}\right]\bu(\bq),\quad 1\leqslant j \leqslant N,
 \end{equation}
 where $\bu(\bq) = \BA(\bq)^{-1}\bff$ is the discrete solution of the boundary value problem (\ref{bvp1}) and $\bv(\bq)$ is the solution to the following adjoint problem
 \begin{equation}
  \BA^*(\bq)\bv(\bq) = -\BC^*\br(\bq),\quad\br(\bq) = \BF(\bq)-\bz
 \end{equation}

The gradient $g_{reg}(q)$ of the ``regularization term'' $\dfrac{1}{2}||\bd^k-\nabla q-\bb^k||_2^2$ is approximated by \cite{vogel1999sparse}
 \begin{equation}\label{greg}
  g_{reg}(q)=\CL(q)q
 \end{equation}
 where $\CL(q)$ is the negative Laplacian operator, i.e. $\CL=-\Delta$, which is in fact independent of $q$. The gradient $\bg_{reg}(\bq)$ of the discretized ``regularization term'' is then obtained via the discretization of (\ref{greg}).
 
 The gradient of the original functional is thus the sum:
 \begin{equation}
  \bg(\bq)=\mu \bg_{ls}(\bq)+\lambda \bg_{reg}(\bq).
 \end{equation}

\subsection{Gradient field estimation on triangle meshes}

The generalized shrinkage formula (\ref{shrinkage}) in the split Bregman scheme requires estimations of gradient fields (with respect to $x$) $\nabla q$. When FEM is applied, the estimations are on the triangular mesh $\CT$, evaluated at mesh points $\bp\in\CP$. A two-step scheme \cite{mancinelli2018gradient} is applied to obtain these estimations.
 
 \medskip
 
 \noindent Step 1. Per-cell linear estimation (PCE)
 
 \medskip
 
 On each triangle $T\in\CT$ with barycentric coordinates \cite{atkinson2005theoretical} \cite{strang2007computational} $\{\varphi^j\}_{j=1}^{3}$ we have
 \begin{equation}
  q|_T=\sum_{j=1}^3(Q_e)_j\varphi^j
 \end{equation}
 where $\BQ_e=((Q_e)_j)_{j=1}^3$ contains the components of $\BQ$ evaluated at the three vertices $\{x_j\}_{j=1}^3$. Then it is natural to consider
 \begin{equation}
  \label{pce}\nabla q|_T=\sum_{j=1}^3(Q_e)_j\nabla\varphi^j.
  \end{equation} 
  (\ref{pce}) can be evaluated without explicitly solving for the barycentric coordinates:
 \begin{equation}
  \nabla q|_T=((Q_e)_2-(Q_e)_1)\frac{(x_1-x_3)^\perp}{2Area(T)}+((Q_e)_3-(Q_e)_1)\frac{(x_2-x_1)^\perp}{2Area(T)}
 \end{equation}
 where $^\perp$ denotes rotating a vector by $\frac{\pi}{2}$.
 
 \medskip
 
 \noindent Step 2. Average gradient on star (AGS)
 
 \medskip
 
 The star $\CN(p)$ of a mesh point $p\in\BP$ is defined as the collection of all triangles $T\in\CT$ that are incident at $p$, that is, all the triangles sharing the common vertex $p$. Then
 \begin{equation}
  \nabla q(p)\approx\frac{1}{\sum_{T\in\CN(p)}Area(T)}\sum_{T\in\CN(p)}Area(T)\nabla q|_T.
 \end{equation}
 
 \subsection{Regularization parameter selection: the L-curve method}
 
 In the split Bregman scheme (Algorithm \ref{alg_1}), there are two regularization parameters to be selected: $\mu$, or equivalently, $\alpha$; and $\lambda$. Empirically, the selection of $\lambda$ depends on the selection of $\alpha$, hence the proposed procedure of regularization parameter selection would be:
 \begin{enumerate}
  \item For each fixed $\alpha$, select $\lambda=\lambda(\alpha)$ as the optimal selection for $\lambda$ under $\alpha$;
  \item Select the optimal $\alpha$ with the $\lambda(\alpha)$ selected previously.
 \end{enumerate}
 In practice, however, this process would require a large amount of computational resources. In this paper, we will only perform the selection of $\lambda$ for fixed $\alpha$.
 
 For each fixed $\alpha$, if $\lambda$ is too small, then the $q$-subproblem (\ref{state1}), which is a penalized least-squares problem in nature, is under-regularized, and the reconstructed $q$ may be highly oscillatory; whereas if $\lambda$ is too large then the reconstructed $q$ may be too smooth and hence inaccurate.
 
 This phenomenon is shown by the L-curve \cite{vogel2002computational}, in which the log of the residual $||\CF(q)-z||_2^2$ is plotted against the log of $||\nabla q||_2^2$ for a range of values of $\lambda$. This curve typically has an L shape. Then the criterion for the selection of $\lambda$ is to pick the value of $\lambda$ corresponding to the ``corner'' of this curve.
 
 In practice, the computation time consumption for L-curve can be reduced using parallel computing in MATLAB or GNU Octave.

 \section{Results}
 
 We present some examples of 2D conductivity reconstruction given by the boundary value problem (\ref{bvp1}) with different geometries of the domain, different boundary conditions, and different forcing functions to demonstrate the effectiveness of the split Bregman scheme. For these examples, simulated data is generated by first solving the boundary value problem using CCFD or FEM accordingly and applying the discrete version of the observation operator $\CC$, then adding Gaussian error with mean zero and with standard deviation selected so that the noise-to-signal ratio is 0.01 \cite{vogel1999sparse}.
 
 All numerical computations were done on a laptop computer with 2.6 GHz up to 4.5 GHz, six cores Intel x86 CPU, and 16GB RAM. 
 
 \subsection{A simple example: L-curve, reconstruction quality, and convergence}

 In this subsection, we illustrate the L-curve method for choosing the parameter $\lambda$, the rate of convergence of the split Bregman scheme, and the quality of reconstruction under different choices of $\lambda$ by a simple example similar to what \cite{vogel1999sparse} has used.
 
 Consider the square domain $\Omega=[0,1]\times[0,1]$ with Dirichlet boundary condition $u=0$ on the left ($x=0$) and right ($x=1$) edges, and no flux Neumann boundary condition $\frac{\partial u}{\partial y}=0$ on the top ($y=1$) and bottom ($y=0$) edges. The forcing function $f$ is given by a point source at $(x_0,y_0)=(\frac{1}{2},\frac{3}{5})$, that is, $f(x,y)=\delta(x-x_0,y-y_0)$.
 
 The conductivity $\kappa(x,y)$ is given by a piece-wise constant function:
 \begin{equation}
  \kappa(x,y) = \left\{
  \begin{array}{ll}
   \kappa_0, & y < \frac{1}{2}, \\
   \kappa_1, & y \geqslant \frac{1}{2}.
  \end{array}
 \right.
 \end{equation}
 In this example we choose $\kappa_0=1$ and $\kappa_1=0.1$.

 Uniform CCFD discretization is applied to the square domain $\Omega$ with $n_x=n_y=50$ cells on each side, so that the total number of cells is $N=2500$.
 
 \subsubsection{L-curve and reconstruction qualities under different choices of parameters}

 \begin{figure}[t]
  \centerline{\includegraphics[width=.55\textwidth]{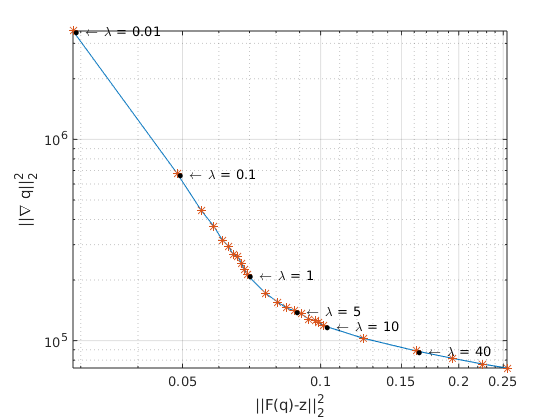}}
  \caption{L-curve for $\alpha=0.0002$, which is a log-log plot of the residual against the square of $L^2$ norm of reconstructed $q$ for different values of $\lambda$.}
  \label{lcurve}
 \end{figure}

 Figure \ref{lcurve} shows the L-curve for choosing the parameter $\lambda$ with fixed $\alpha=0.0002$. This L-curve yields $\lambda=5$ as the choice of the parameter $\lambda$.

\begin{figure}[t]
  \centerline{\includegraphics[width = \textwidth]{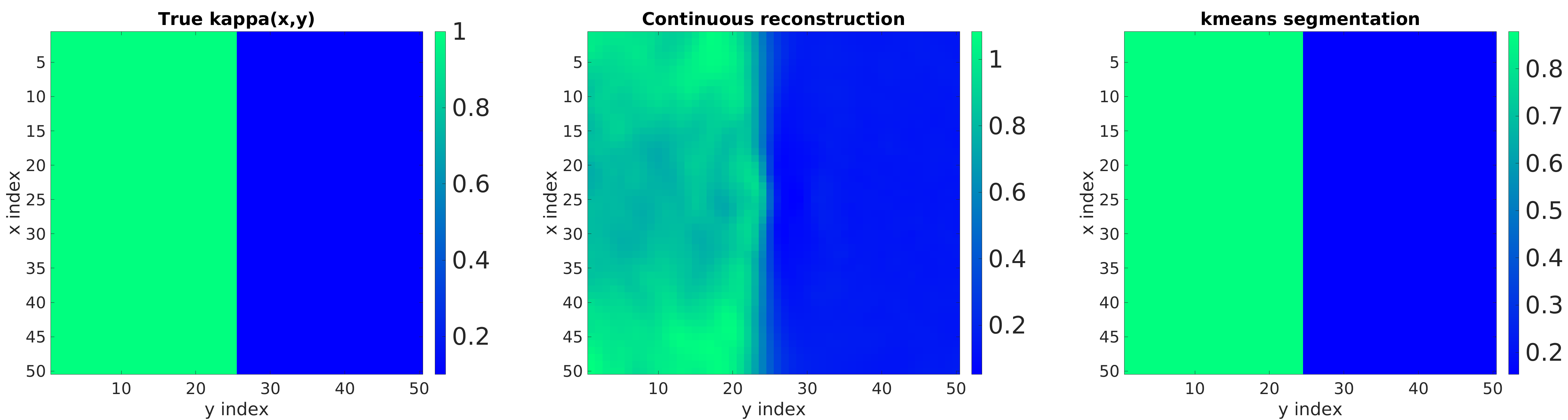}}
  \caption{Reconstruction results for $\lambda=5$: true conductivity (left), continuous reconstruction (middle), and \texttt{kmeans} segmentation (right).}
  \label{Ex5}
 \end{figure}
\begin{figure}[t]
  \centerline{\includegraphics[width = \textwidth]{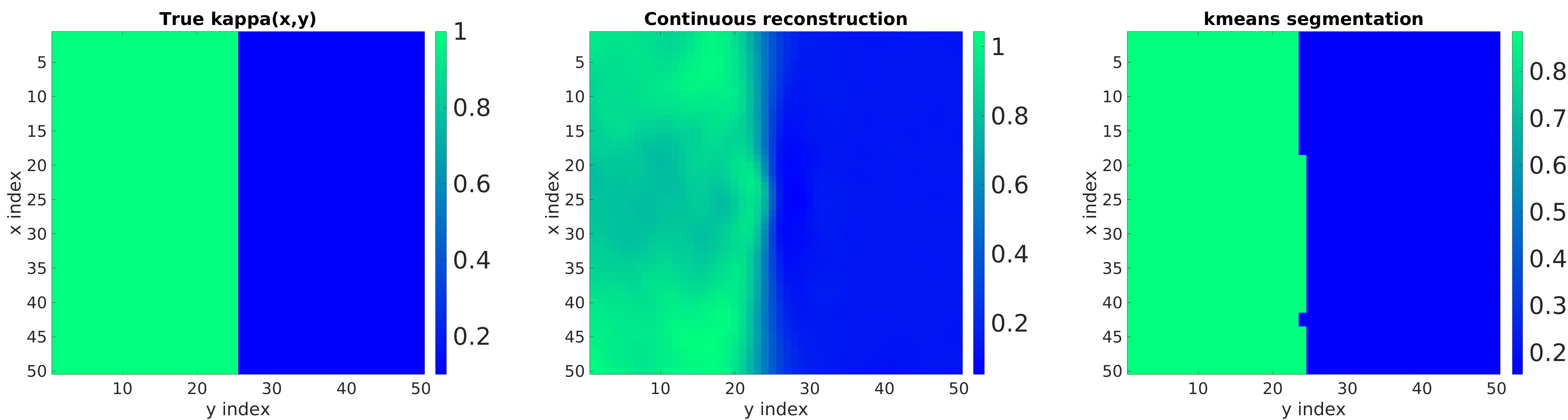}}
  \caption{Reconstruction results for $\lambda=10$: true conductivity (left), continuous reconstruction (middle), and \texttt{kmeans} segmentation (right).}
  \label{Ex10}
 \end{figure}

Figure \ref{Ex5} and Figure \ref{Ex10} show the true conductivity $\kappa$, the continuous reconstruction, and the result of the segmentation by \texttt{kmeans}, with $\lambda=5$ and $\lambda=10$, respectively. Although the continuous reconstruction with $\lambda=10$ seems to be less noisy than that with $\lambda=5$, the quality of the interface reconstruction after \texttt{kmeans} is significantly higher with $\lambda=5$.

 \subsubsection{Convergence of the split Bregman scheme}
 
 We use normalized relative error \cite{goldstein2009split} to specify the termination criterion and show convergence of the split Bregman iterations:
 \begin{equation}
  err^k = \frac{||q^k-q^{k-1}||_2^2}{||q^k||_2^2}.
 \end{equation}
 
 \begin{figure}[t]
  \centerline{\includegraphics[scale=0.5]{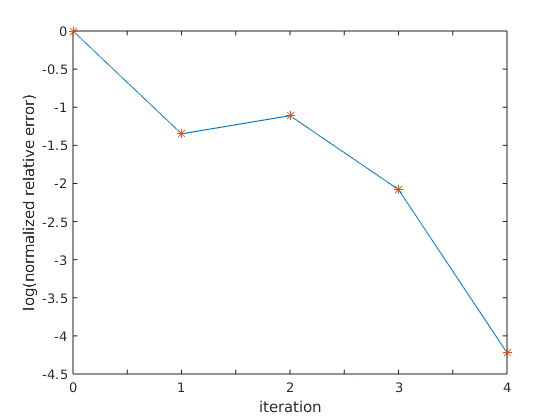}}
  \caption{Normalized relative error vs. iteration number for the split Bregman iterations with $\alpha=0.0002$ and $\lambda=5$.}
  \label{Convergence}
 \end{figure}
 
 Figure \ref{Convergence} shows the rate of convergence of the split Bregman iterations with $\alpha=0.0002$ and $\lambda=5$. The termination criterion is set to be $err^k<10^{-6}$. The iteration terminated at the fifth iteration. 
 
 As is shown by the figure, the split Bregman iteration converges very fast in this example of conductivity reconstruction. Note that semiconvergence behavior is expected if the data is noisy \cite{osher2005iterative}.
 
 \subsection{Simple domain geometry with complicated interface geometry between different phases of conductivity}
 
 In the next example, we demonstrate the performance of the split Bregman scheme applied to conductivity reconstruction problems in which the geometries of the interfaces between different phases of conductivity.
 
 Consider the square domain $\Omega=[-\frac{1}{2},\frac{1}{2}]\times[-\frac{1}{2},\frac{1}{2}]$, with the same boundary conditions and forcing function (which is shifted) as in the previous example. The conductivity is given by a piece-wise constant function
 \begin{equation}
  \kappa(x,y) = \left\{
  \begin{array}{ll}
   \kappa_0, & (x,y)\in\Omega_0, \\
   \kappa_1, & (x,y)\in\Omega\setminus\Omega_0,
  \end{array}
 \right.
 \label{kappa}
 \end{equation}
 where $\kappa_0=1$ and $\kappa_1=0.1$ as before, and $\Omega_0$ is chosen to have the shape of a four-leaved clover, the boundary of which is given in polar coordinates by
 \begin{equation}
  r = \left[\left(\frac{1}{2}\sin 2\theta + \frac{1}{8}\sin 6\theta\right)^4+\epsilon\right]^\frac{1}{4},
  \label{clovereq}
 \end{equation}
 where $\epsilon = 10^{-3}$ is a small number intended to ``soften'' the sharp corners of the four leaves; otherwise the point $(x,y)=(0,0)$ would be a singular point.
 \begin{figure}[t]
  \centerline{\includegraphics[width = \textwidth]{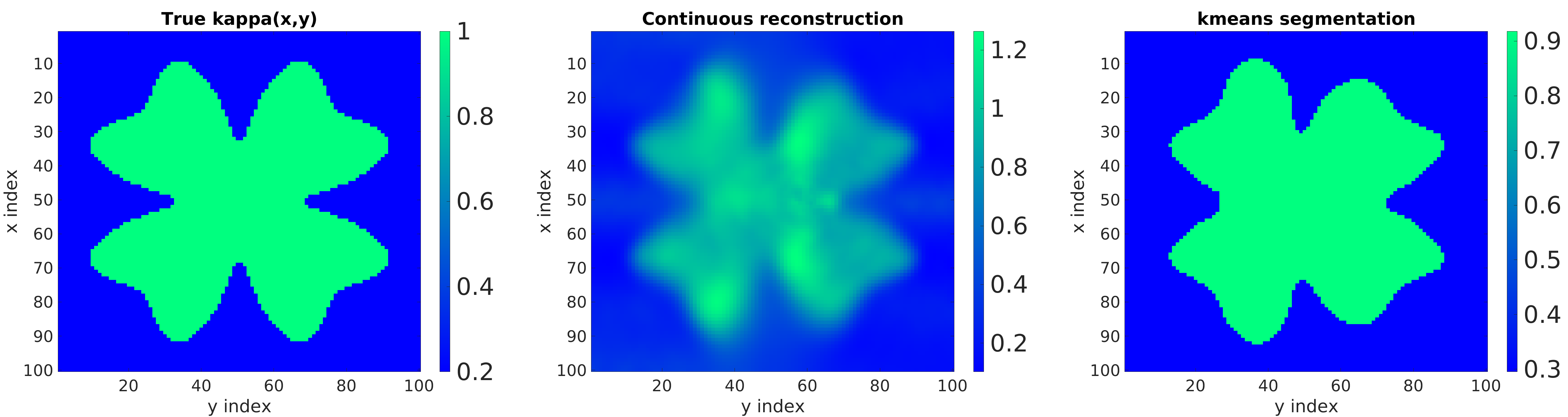}}
  \caption{Reconstruction results for the clover-shaped interface example: true conductivity (left), continuous reconstruction (middle), and \texttt{kmeans} segmentation (right). Number of nodes is $N=10000$; $\alpha=0.0001$, $\lambda=15$.}
  \label{Clover}
 \end{figure}
 
 Figure \ref{Clover} shows the reconstruction result with a uniform CCFD discretization with $n_x=n_y=100$ cells on each side, and with regularization parameters $\alpha=0.0001$ and $\lambda=15$. From the results, we see that the continuous reconstruction obtained directly from the split Bregman scheme is able to recover the profile of the clover-shaped subset $\Omega_0$. 
 
 
 
 \subsection{Complicated domain geometry}
 
 In the following examples, we apply the split Bregman scheme to non-rectangular domains, in which case FEM methods have to be used. Dirichlet boundary conditions are imposed to the boundary value problem (\ref{bvp1}), with forcing function $f=0$.
 
 The conductivity $\kappa$ is given by (\ref{kappa}), with $\kappa_0=1$ and $\kappa_1=0.1$ and $\Omega_0$ specified for each particular example.
 
 \subsubsection{Disc domain}
 
 In this example, the domain $\Omega$ is chosen to be a disc with radius $1$, centered at the origin. Periodic boundary condition given by polar coordinates
 \begin{equation}
  g_D(r,\theta) = r^3\sin(3\theta)
 \end{equation}
 is imposed. The subset $\Omega_0$ is set to be all the nodes $p$ that are in the open disc with radius $\frac{1}{2}$, centered at the origin.
 
 \begin{figure}[t]
  \centerline{\includegraphics[width = \textwidth]{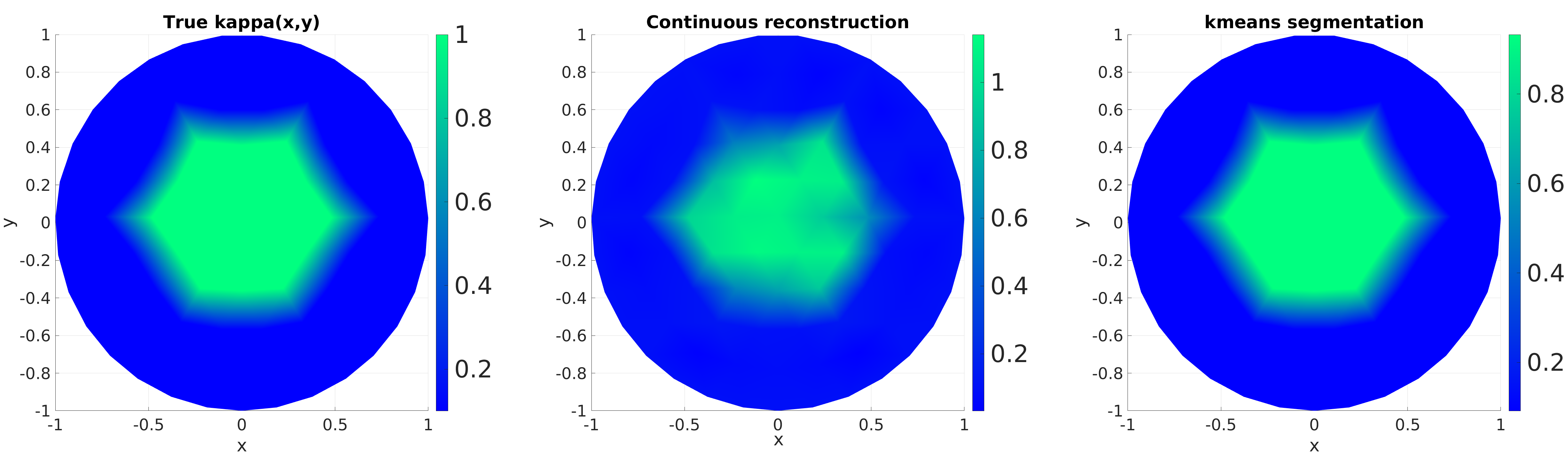}}
  \caption{Reconstruction results for the disc domain example: true conductivity (left), continuous reconstruction (middle), and \texttt{kmeans} segmentation (right). $\alpha=0.0005$, $\lambda=1$.}
  \label{Disc}
 \end{figure}

 Figure \ref{Disc} shows the reconstruction results with parameters $\alpha=0.0005$ and $\lambda=1$.
 
 \subsubsection{Annular domain}
 
 Next, we consider a slightly more complicated domain $\Omega$, which is an annular domain given by removing a disc of radius $0.4$ from a disc of radius $1$. The boundary condition is given by $g_D=0.1$ on the inside boundary and $g_D=0.5$ on the outside boundary. The subset $\Omega_0$ is set to be all the nodes $p$ that lie in the open disc with radius $0.7$, centered at the origin.
 
  \begin{figure}[t]
  \centerline{\includegraphics[width = \textwidth]{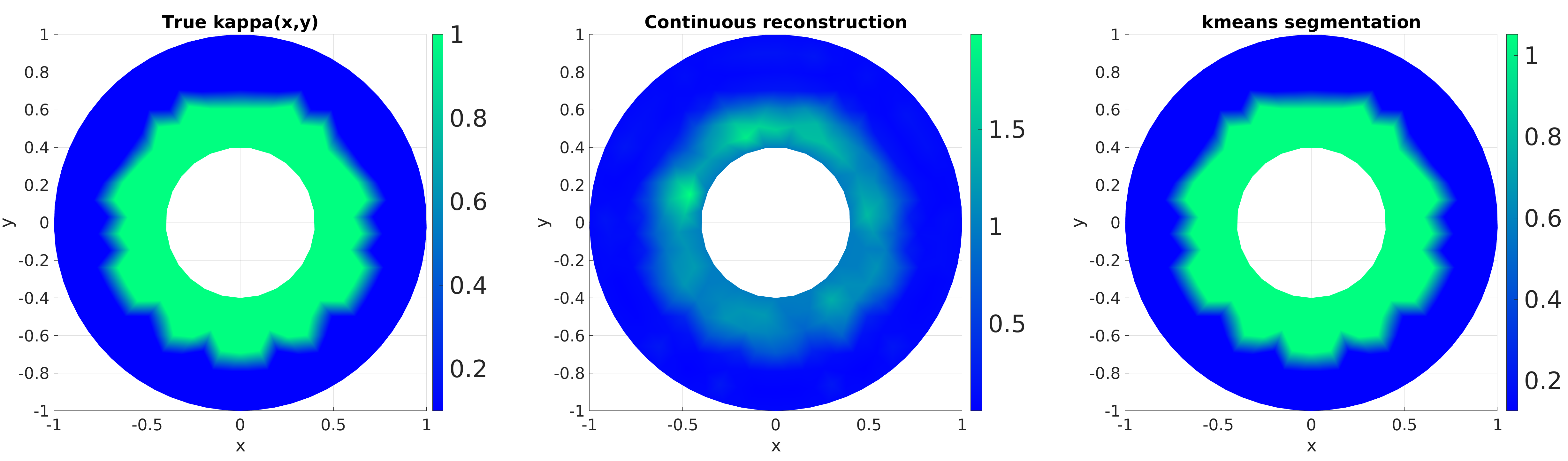}}
  \caption{Reconstruction results for the annular domain example: true conductivity (left), continuous reconstruction (middle), and \texttt{kmeans} segmentation (right). $\alpha=0.0005$, $\lambda=14$.}
  \label{Annulus}
 \end{figure}

 Figure \ref{Annulus} shows the reconstruction results with parameters $\alpha=0.0005$ and $\lambda=14$.
 
 \subsubsection{More complicated domain geometry}
 
 We consider an even more complicated domain $\Omega$ given by the intersection of the following two regions \cite{persson2004simple}:
 \begin{equation}
  y\leqslant \cos x\quad\mbox{and}\quad y\geqslant 5\left(\frac{2x}{5\pi}\right)^4-5.
 \end{equation}
 The boundary condition is given by $g_D=0.1$ on the lower boundary and $g_D=0.9$ on the upper boundary. The subset $\Omega_0$ is set to be all the nodes $p$ that lie in the open disc with radius $1.4$, centered at the $(0,-2.5)$.
 
 \begin{figure}[t]
  \centerline{\includegraphics[width = \textwidth]{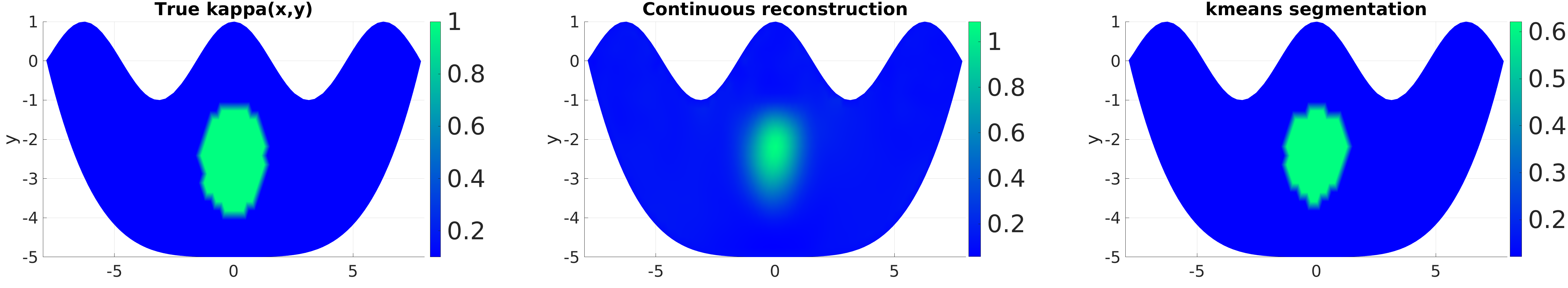}}
  \caption{Reconstruction results for the ``crown-shaped'' domain example: true conductivity (left), continuous reconstruction (middle), and \texttt{kmeans} segmentation (right). $\alpha=0.0001$, $\lambda=150$.}
  \label{Crown}
 \end{figure}
 Figure \ref{Crown} shows the reconstruction results with parameters $\alpha=0.0005$ and $\lambda=14$.
 
 \section{Discussion}
 
 In the previous section, we have shown that the split Bregman scheme can be applied to conductivity reconstruction problems with various domain geometries and get good reconstruction results. Working together with a \texttt{kmeans} segmentation process, the scheme can further recover the interfaces between different phases of ``blocky'' conductivity functions. Still, there are several aspects in which the current reconstruction scheme can be improved, and this reconstruction scheme may also be applied to other types of parameter reconstruction problems.
 
 \subsection{Improving the numerical optimization efficiency and accuracy}
 
 In the $q$-subproblem \ref{state1} of the split Bregman iteration scheme (Algorithm \ref{alg_1}), there is no $L^1$ portion, so it is differentiable. However, it is inherently a Tikhonov-regularized nonlinear least-squares minimization problem that has the same structure as (\ref{pls1}). This means that the $q$-subproblem is still a quite challenging optimization problem and is difficult to solve efficiently and accurately using BFGS-type methods, especially when the spatial grid is fine and the number $N$ of discrete nodes is large. 
 
 The robustness of the split Bregman scheme guarantees that the iterations will converge even if the solution to the $q$-subproblem is only approximate. Nevertheless, a better optimization solver for the $q$-subproblem is still desired to further improve the quality of the reconstructions and to make the split Bregman scheme more efficient.
 
 Alternatively, it may be possible to reduce the number of variables in the numerical optimization solver, while keeping the spatial grid fine with large $N$. To achieve this, we would need to apply model order reduction techniques (see, for example, \cite{schilders2008model}) to the log conductivity function $q(x)$, but doing so may also require modifications of the split Bregman iteration scheme.
 
 \subsection{Improving the segmentation scheme}
 
 The clustering and thresholding scheme for ``blocky'' conductivity functions may be improved as well. In the current scheme, fine tuning of the \texttt{kmeans} algorithm is needed sometimes, in order to achieve segmentations with higher precision. So, more accurate clustering algorithms may be used in the future to further improve the precision of the segmentations.
 
 As we mentioned before, the current segmentation scheme cannot automatically select the number $K$ of different clusters. An automatic selection algorithm for $K$ would need to analyze the structure of the reconstructed $\bkappa(p)$ and recognize its pattern, which would lead us to the realm of pattern recognition.
 
 \subsection{Better methods for regularization parameter selection}
 
 To the best of our knowledge, currently there is no automatic regularization parameter selection method, designed specifically for the split Bregman iteration scheme, that can simultaneously select the optimal values for $\mu$ and $\lambda$. 
 
 That being said, there exists general parameter selection methods that can simultaneously select more than one regularization parameters. One example is the bilevel optimization method \cite{antil2020bilevel}. Briefly speaking, if we have a collection of true log conductivity function samples $\left\{\bq_{true}^{(i)}\right\}_{i=1}^{N_s}$, we can solve the following minimization problem
 \begin{equation}\label{bilevel}
  \min_{\mu,\lambda}\frac{1}{2N_s}\left\|\bq^{(i)}(\mu,\lambda)-\bq_{true}^{(i)}\right\|_2^2,
 \end{equation}
 where $\bq^{(i)}(\mu,\lambda)$ is the reconstruction result corresponding to the sample $\bq_{true}^{(i)}$. Clearly, (\ref{bilevel}) is nonlinear, and it would be difficult to apply gradient based optimization methods to (\ref{bilevel}).
 
 On the other hand, it is natural to ask whether the reconstruction results obtained from the split Bregman scheme are sensitive to the regularization parameters. Although our results show that changing $\mu$ and $\lambda$ will indeed affect the reconstruction results, to the best of our knowledge, currently there is no theoretical analysis. Thus it may be of interest to take a deeper look at the parameters in the split Bregman method.
 
 \subsection{Generalization to other types of parameter reconstruction problems}
 
 In this paper, we have only considered the reconstruction of the conductivity functions in steady-state, elliptic diffusion equations. It may be possible to apply this split Bregman scheme to parameter reconstruction problems for other types of equations. For example, we could consider the time-dependent, parabolic diffusion equations corresponding to (\ref{bvp1}). 
\printbibliography

\appendix

\end{document}